\documentclass[12pt]{article}

\usepackage{lineno,hyperref}
\usepackage{amssymb}
\usepackage{amsmath}
\usepackage{amsthm}
\usepackage{latexsym}
\usepackage{amsfonts}
\usepackage{graphicx}
\usepackage{graphics}

\theoremstyle{plain}
 \newtheorem{theorem}{Theorem}
 \newtheorem{lemma}[theorem]{Lemma}
  \newtheorem{proposition}[theorem]{Proposition}
 
  \newtheorem{remark}[theorem]{Remark}

\theoremstyle{definition}
\newtheorem*{Proof}{Proof}
\newcommand{\fa}{\forall}
\newcommand{\dis}{\displaystyle}

\newcommand{\ca}{\mathcal{A}}

\newcommand{\ce}{\mathcal{E}}

\newcommand{\cf}{\mathcal{F}}
\newcommand{\cg}{\mathcal{G}}

\newcommand{\ld}{\ldots}

\newcommand{\el}{\ell}

\newcommand{\ra}{\rightarrow}

\newcommand{\al}{\alpha}
\newcommand{\bi}{\beta}
\newcommand{\ga}{\gamma }

\newcommand{\de}{\delta }

\newcommand{\e}{\varepsilon }

\newcommand{\thi}{\theta }

\newcommand{\mi}{\mu }

\newcommand{\R}{\mathbb{R}}

\newcommand{\N}{\mathbb{N}}

\newcommand{\ssum}{\sum\limits}

\newcommand{\qs}{$\quad\square$}
\newcommand{\hi}{\widehat{I}}

\newcommand{\dil}{\mbox{\em dil}}

\begin{document}
\title{On the speed of convergence in the strong\\ density theorem}
\author{Panagiotis Georgopoulos, Constantinos Gryllakis$^\ast$
\\
Department of Mathematics,\\ National and Kapodistrian University of Athens,\\
Panepistimiopolis 15784, Athens, Greece}\footnotetext{\hspace*{-0.5cm}$^\ast$
Corresponding author\\
{\em E-mail addresses}: pangeorgopoul@gmail.com \& pgeorgop@math.uoa.gr
(P. Georgopoulos), cgryllakis@math.uoa.gr (C. Gryllakis)}
\date{}
\maketitle
\begin{abstract}
For a compact set $K\subseteq\R^m$, we have two indexes given under simple parameters of the set $K$ (these parameters go back to Besicovitch and Taylor in the late 50's).
In the present paper we prove that with the exception of a single extreme value for each index, we have the following elementary estimate on how fast the ratio in the strong density theorem of Saks will tend to one
\[
\frac{|R\cap K|}{|R|}>1-o\bigg(\frac{1}{|\log d(R)|}\bigg) \qquad \text{for a.e.} \ \ x\in K \ \ \text{and for} \ \ d(R)\ra0
\]
(provided $x\in R$, where $R$ is an interval in $\R^m$, $d$ stands for the diameter and $|\cdot|$ is the Lebesgue measure).\\
This work is a natural sequence of \cite{3} and constitutes a
contribution to Problem 146 of Ulam [5, p. 245] (see also [8, p.78])
and Erd\"{o}s' Scottish Book `Problems' [5, Chapter 4, pp. 27-33],
 since it is known that no general statement can be made on how fast
the density will tend to one.
\end{abstract}
\section{Introduction}\label{sec1}
\noindent

Given a compact set $K\subseteq\R^2$, it has always the form $K=I\big\backslash\dis\bigcup_{n\in\N}(I_n\times J_n)$, where $I=(\al,b)\times(c,d)$ and $I_n\times J_n$, $n\in\N$ are disjoint cubes $(|I_n|=|J_n|)$, with $|I_n|=w_n$, $n\in \N$, in a non-increasing order. For the sequence $\{w_n:
\linebreak n\in\N\}$, we write
\[
r_n=\sum^\infty_{m=n}w^2_m, \ \ n\bigg(\frac{r_n}{n}\bigg)^{a_n}=1, \ \ a\{w_n^2:n\in\N\}=\underset{n}{\lim\inf}a_n.
\]
Clearly $0\le a\{w^2_n:n\in\N\}\le1$.\\

Second, we have the well-known Besicovitch-Taylor index (or exponent of convergence) of the sequence $\{w^2_n:n\in\N\}$
\[
e_{BT}\{w^2_n:n\in\N\}:=\inf\{c>0:\sum^\infty_{n=1}(w^2_n)^c\;\text{converges}\}
\]
(see [7, p. 34 and p. 292]).\\
Clearly $0\le e_{BT}\{w^2_n:n\in\N\}\le1$.

Note that both of the above indexes used by Besicovitch and Taylor, in order to estimate the Hausdorff dimension of some exceptional set (see [1,2]).

In the present paper, by an essential modification of the methods of \cite{3}, we prove that for a compact set $K$ of positive Lebesgue measure with
\[
a\{w^2_n:n\in\N\}\neq 0 \ \ \text{and} \ \ e_{BT}\{w^2_n:n\in\N\}\neq1
\]
we have the following elementary estimate, on how fast the density tends to one
\[
\inf_{x\in A\times B\atop diam(A\times B)<t}\frac{|(A\times B)\cap K|}{|A\times B|}\ge1-o\bigg(\frac{1}{|\log t|}\bigg)
\]
for a.e. $x\in K$ and for sufficiently small $t$.

Note that our methods work equally well in any dimension, but for simplicity reasons
we restrict to $\R^2$ and it is worth mentioning that the estimate is independent of the dimension.

This is a contribution to Problem 146 of Ulam [5, p. 245] (see also [8, p. 78]) and Erd\"{o}s' Scottish Book `Problems' [5, Chapter 4, pp. 27-33]. It is known that no general statement can be made on how fast the density tends to one, although it is known that for a particular set there exists a function, depending on the set, that is dominated by the density (see \cite{3,6}). However, it is not clear how this particular function
depends upon the set (note that the proof in \cite{6} is non-constructive).

\section{The notion of simultaneous dilation}\label{sec2}
\noindent

We recall the notion of simultaneous dilation from \cite{3} (see Sections 2 and 3) and assume all the curriculum of propositions and lemmata from there.

Let $\{I_i:i\in\mathcal{A}\}$ be a finite pairwise disjoint collection of bounded intervals in $\R$ and let $\ga>1$.\\
We enumerate this collection as $I_k$, $k=1,\ld,n$, such that
\[
\inf I_1\le\sup I_1\le\inf I_2\le\sup I_2\le\cdots\le\inf I_{n-1}\le\sup I_{n-1}\le\inf I_n\le\sup I_n.
\]

We define inductively for $k=1$, $I'_1$ an open interval with the same right end as $I_1$ such that
\[
\big|I'_1\big\backslash\bigcup^n_{i=1}I_i\big|=\ga\cdot|I_1|
\]
and an open interval $I''_1$ with the same left end as $I_1$ such that
\[
\big|I''_1\big\backslash\bigcup^n_{i=1}I_i\big|=\ga\cdot|I_1|.
\]
We set $\widehat{I}_1:=I'_1\cup I''_1$.

If $\hi_1,\ld,\hi_k$, $k<n$ are defined, we define $I'_{k+1}$ an open interval with the same right end as $I_{k+1}$ such that
\[
\big|I'_{k+1}\big\backslash\big(\bigcup^k_{i=1}\hi_i\cup\bigcup^n_{i=k+1}I_i\big)\big|=\ga\cdot|
I_{k+1}|
\]
and $I''_{k+1}$ an open interval with the same left end as $I_{k+1}$ such that
\[
\big|I''_{k+1}\big\backslash\big(\bigcup^k_{i=1}\hi_i\cup\bigcup^n_{i=k+1}I_i\big)\big|=\ga\cdot
|I_{k+1}|.
\]
We set $\hi_{k+1}:=I'_{k+1}\cup I''_{k+1}$.\\
We shall call $\bigcup\limits^n_{k=1}\hi_k$ the {\em one-dimensional simultaneous $\ga$-dilation} of the union $\bigcup\limits_{i\in\ca}I_i$, or simply the $\ga$-{\em dilation}. Symbolically
\[
\bigcup^n_{k=1}\hi_k:=\ga-\dil\big(\bigcup_{i\in\ca}I_i\big).
\]

Next, concerning the two-dimensional dilation, we consider $\{I_i\!\times\! J_i:\;i\!\in\!\ca\}$ a finite pairwise disjoint collection of bounded intervals in $\R^2$ and a $\ga>1$.

For every non-empty $\ce\subset\ca$ we set
\[
\begin{array}{l}
  R_\ce:=\{x\in\R:\;x\in I_i\Leftrightarrow i\in\ce, \ \ \fa\;i\in\ca\} \\ [1ex]
  Q_\ce:=\{y\in\R:\;y\in J_i\Leftrightarrow i\in\ce, \ \ \fa\;i\in\ca\}.
\end{array}
\]

It is clear that there exists an $\cf\subseteq2^\ca\times2^\ca$ ($2^\ca$ denotes the power set of $\ca$) such that for $(\al,\bi)\in\cf$, $R_\al\times Q_\bi\neq\emptyset$ and
\[
\bigcup_{i\in\ca}(I_i\times J_i)=\bigcup_{(\al,\bi)\in\cf}(R_\al\times Q_\bi),
\]
(note that the respective sides of $R_\al\times Q_\bi$' s are pairwise disjoint).\\
We set for $\bi\in2^\ca$
\[
\cf_\bi:=\{\al\in2^\ca:\;(\al,\bi)\in\cf\}, \quad \cf':=\{\bi\in2^\ca:\cf_\bi\neq\emptyset\}
\]
and for $\bi\in\cf'$
\[
D_\bi:=\ga-\dil\big(\bigcup_{\al\in\cf_\bi}R_\al\big)
\]
(the one dimensional simultaneous $\ga$-dilation of $\bigcup\limits_{\al\in\cf_\bi}R_\al$).

For every non-empty $\cg\subseteq\cf'$ we set
\[
V_\cg:=\{x\in\R:\;x\in D_\bi\Leftrightarrow\bi\in\cg, \ \ \fa\;\bi\in\cf'\}
\]
and
\[
H_\cg:=\ga-\dil\big(\bigcup_{\bi\in\cg}Q_\bi\big).
\]

We shall call $\bigcup\limits_{\cg\subseteq\cf'}(V_\cg\times H_\cg)$ the {\em two dimensional simultaneous $\ga$-dilation}, or simply the $\ga$-{\em dilation} of the union $\bigcup\limits_{i\in\ca}(I_i\times J_i)$. Symbolically
\[
\bigcup_{\cg\subseteq\cf'}(V_\cg\times H_\cg):=\ga-\dil\big(\bigcup_{i\in\ca}(I_i\times J_i)\big).
\]
\section{The auxiliary function}\label{sec3}
\noindent

Throughout this paper, $K$ is a compact subset of $\R^2$ of positive Lebesgue measure, as in the introduction (Section 1).\\
Note that $a\{w^2_n:n\in\N\}\neq0$ and $e_{BT}\{w^2_n:n\in\N\}\neq1$.

It is known (see [4, pp. 274-275] and [7, p. 292 \& THEOREM, p. 35]) that
\[
e_{BT}\{w^2_n:n\in\N\}=\underset{n}{\lim\sup}\frac{\log n}{|\log w^2_n|}.
\]
Since $e_{BT}\{w^2_n:n\in\N\}<1$, for $e_{BT}\{w^2_n:n\in\N\}<\thi<1$ we have finally for every $n\in\N$
\begin{align}
&\frac{\log n}{|\log w^2_n|}<\thi<1\Leftrightarrow\log n<\thi\cdot\log\frac{1}{w^2_n} \nonumber \\
&\Leftrightarrow n<\bigg(\frac{1}{w^2_n}\bigg)^\thi\Leftrightarrow n^{\frac{1}{\thi}}<\frac{1}{w^2_n} \nonumber \\
&\Leftrightarrow w^2_n<\bigg(\frac{1}{n}\bigg)^{\frac{1}{\thi}}.  \label{eq1}
\end{align}

Next, recalling the Bouligand-Minkowski index of the sequence $\{w^2_n:\linebreak n\in\N\}$ (see [7, p. 35])
\[
e_{BM}\{w^2_n:n\in\N\}=\inf\{a:(w^2_n)^{a-1}\cdot\sum^\infty_{i=n}w^2_i \ \  \text{tends to} \ \  0\}
\]
and the fact that this index equals to the Besicovitch-Taylor index $e_{BT}\{w^2_n:n\in\N\}$ (see [7, THEOREM, p. 35]), we have that for
$e_{BT}\{w^2_n:n\in\N\}=e_{BM}\{w^2_n:n\in\N\}<\de<1$
\begin{eqnarray}
r_n\bigg(=\sum^\infty_{m=n}w^2_m\bigg)<(w^2_n)^{1-\de} \qquad \text{finally for every} \ \ n\in\N. \label{eq2}
\end{eqnarray}
So, since $1-\de>0$, by (\ref{eq1}) and (\ref{eq2}) we get
\[
r_n<\bigg(\frac{1}{n}\bigg)^{\frac{1}{\thi}\cdot(1-\de)} \qquad \text{finally for every} \ \ n\in\N
\]
and setting $\dfrac{1}{\thi}\cdot(1-\de)=:\e>0$, we deduce that
\begin{eqnarray}
r_n<\bigg(\frac{1}{n}\bigg)^\e \qquad \text{finally for every} \ \ n\in\N.  \label{eq3}
\end{eqnarray}

Next, we introduce the following sequence of positive integers
\[
\{n_s:=s^s:s\in\N\}
\]
that plays a fundamental role in the definition of the auxiliary function.

We have
\begin{lemma}\label{lem1}
The series $\ssum^\infty_{s=1}2^s\cdot r_{n_s}$ converges.
\end{lemma}
\begin{Proof}
By (\ref{eq3}), it suffices to prove that the series
\[
\sum^\infty_{s=1}2^s\cdot\bigg(\frac{1}{n_s}\bigg)^\e
\]
converges.\\
We shall use the ratio test.\\
We have
\[
\frac{(s^s)^\e}{\big((s+1)^{s+1}\big)^\e}=\bigg(\frac{s^s}{(s+1)^{s+1}}\bigg)^\e=
\bigg(\bigg(\frac{s}{s+1}\bigg)^s\cdot\frac{1}{s+1}\bigg)^\e\ra\bigg(\frac{1}{e}\bigg)^\e\cdot0=0 \]
for $s\ra\infty$.\\
Clearly
\[
\bigg(\frac{n_s}{n_{s+1}}\bigg)^\e=\frac{(s^s)^\e}{\big((s+1)^{s+1}\big)^\e}\ra0 \qquad \text{for} \ \ s\ra\infty.
\]
An immediate application of the ratio test gives the convergence of the\linebreak series. \qs
\end{Proof}

Also, a similar application of the ratio test gives
\begin{lemma}\label{lem2}
The series $\ssum^\infty_{s=1}2^s\cdot\sqrt{r_{n_s}}$ converges.
\end{lemma}
\begin{remark}\label{rem3}
By Lemma \ref{lem2}
\[
2^s\cdot\sqrt{r_{n_s}}\ra0 \qquad \text{for} \ \ s\ra\infty
\]
and since clearly
\[
\sqrt{r_{n_s}}\ge w_{n_k}  \qquad \text{for} \ \ k\ge s
\]
we have
\begin{eqnarray}
2^s\cdot w_{n_{s+1}-1}\ra0  \qquad \text{for} \ \ s\ra\infty.  \label{eq4}
\end{eqnarray}
\end{remark}

Next, we choose a subsequence $\{n_{s_\el}:\el\in\N\}$ of $\{n_s:=s^s:s\in\N]$ as follows
\begin{eqnarray}
2^{s_{\el+1}}\cdot w_{n_{s_{\el+1}+1}-1}\lvertneqq 2^{s_\el}\cdot w_{n_{s_\el+1}-1}  \label{eq5}
\end{eqnarray}
and
\begin{eqnarray}
2^{s_\el}\cdot w_{n_{s_\el+1}-1}\le 2^s\cdot w_{n_{s+1}-1}  \qquad \text{for} \ \ s_\el\le s<s_{\el+1} \label{eq6}
\end{eqnarray}
(i.e. $s_{\el+1}$ is the first positive integer after $s_\el$ such that $2^{s_{\el+1}}\cdot w_{n_{s_{\el+1}+1}-1}\lvertneqq 2^{s_\el}\cdot w_{n_{s_\el+1}-1}$, for $\el\in\N$).

We are now in position to define the auxiliary function $h$ {\em associated to the decomposition of} $K=I\big\backslash\dis\bigcup_{n\in\N}(I_n\times J_n)$ into disjoint cubes, as follows
\[
h(t):=\left\{\begin{array}{l}
               1-\ssum^\infty_{k=1}\dfrac{2}{2^k} \ \  \text{for} \ \ t\ge2^{s_1}\cdot w_{n_{s_1+1}-1}  \vspace*{0.3cm} \\
               1-\ssum^\infty_{k=s_{\el+1}}\dfrac{2}{2^k} \ \ \text{for} \ \ 2^{s_{\el+1}}\cdot w_{n_{s_{\el+1}+1}-1}\le t<2^{s_\el}\cdot w_{n_{s_\el+1}-1}\\
               \hspace*{7cm}  \text{and} \ \ \el\in\N.
             \end{array}\right.
\]
\begin{remark}\label{rem4}
Concerning $h$, we have by (\ref{eq5}) and (\ref{eq6}) that it is well-defined. Also, by (\ref{eq4}), we have $h:(0,+\infty)\ra\R$ and clearly
$\dis\lim_{t\ra0^+}h(t)=1$.
\end{remark}

We shall need the following elementary lemma
\begin{lemma}\label{lem5}
Under the above considerations, the following series
\[
\sum^\infty_{s=1}2^s\cdot r_{n_s}, \ \ \sum^\infty_{s=1}2^s\cdot\sqrt{r_{n_s}} \ \ \text{and} \ \ \sum^\infty_{s=1}(2^s)^2\cdot r_{n_s}
\]
converge.
\end{lemma}
\begin{Proof}
The convergence of the first two series is already given by Lemmata \ref{lem1} and \ref{lem2}. The convergence of the third one also follows easily from the ratio test. \qs
\end{Proof}

Under the above considerations we have
\begin{theorem}\label{thm6}
Let $K\subseteq\R^2$ be a compact set of positive Lebesgue measure, $I=(\al,b)\times(c,d)$ such that $K\subseteq I$ and $K=I\big\backslash\dis\bigcup_{n\in\N}(I_n\times J_n)$, where $I_n\times J_n$, $n\in\N$ are disjoint cubes, $I_n=J_n=w_n$, with $\{w_n:n\in\N\}$ a non-increasing sequence. Then, for $h$ associated to the decomposition of $K=I\big\backslash\dis\bigcup_{n\in\N}(I_n\times J_n)$, we have
\[
\inf_{x\in A\times B\atop diam(A\times B)<t}\frac{|(A\times B)\cap K|}{|A\times B|}\ge h(t)
\]
for every\;\;
$x\in K\big\backslash\big(\big(\bigcup\limits^\infty_{n=1}\partial(I_n\times J_n\big)\big)\cup\big(\bigcap\limits^\infty_{m=1}\bigcup\limits^\infty_{s=m}2^s-dil
\big(\bigcup\limits^{n_{s+1}-1}_{i=n_s}(I_i\times J_i)\big)\big)\big)$\\
and for sufficiently small $t<\de(x)$, where $A\times B$ is a bounded interval in $\R^2$ ($\partial$ denotes the topological boundary).

Moreover
\[
\big|\big(\bigcup^\infty_{n=1}\partial(I_n\times J_n\big)\big)\cup\big(\bigcap^\infty_{m=1}\bigcup^\infty_{s=m}2^s-dil\big(
\bigcup^{n_{s+1}-1}_{i=n_s}(I_i\times J_i)\big)\big)\big|=0.
\]
\end{theorem}
\begin{Proof}
Firstly, we prove the following \vspace*{0.2cm} \\
\noindent
{\bf Claim}. $\big|\big(\bigcup\limits^\infty_{n=1}\partial(I_n\times J_n\big)\big)\cup\big(\bigcap\limits^\infty_{m=1}\bigcup\limits^\infty_{s=m}2^s-dil\big(
\bigcup\limits^{n_{s+1}-1}_{i=n_s}(I_i\times J_i)\big)\big)\big|=0$.\bigskip

For brevity in the notation we set for every $m\in N$
\[
C_m:=\bigcup^\infty_{s=m}2^s-dil\big(\bigcup^{n_{s+1}-1}_{i=n_s}(I_i\times J_i\big)\big).
\]
By [3, Proposition 3.2]
\[
\bigg|2^s-dil\bigg(\bigcup^{n_{s+1}-1}_{i=n_s}(I_i\times J_i)\bigg)\bigg|=(2\cdot2^s+1)^2\cdot\sum^{n_{s+1}-1}_{i=n_s}w^2_i\le(2\cdot
2^s+1)^2\cdot r_{n_s}
\]
so
\[
|C_m|\le\sum^\infty_{s=m}(2\cdot2^s+1)^2\cdot r_{n_s} \qquad \text{for} \ \ m\in\N.
\]
Since by Lemma \ref{lem5}, the series $\ssum^\infty_{s=1}(2\cdot 2^s+1)^2\cdot r_{n_s}$ converges, we have
\[
\big|\bigcap^\infty_{m=1}C_m\big|=0.
\]
Clearly
\[
\big|\bigcup^\infty_{n=1}\partial(I_n\times J_n)\big|=0,
\]
so the claim holds true.\bigskip

Next, we take some $x\in K$ so that
\[
x\notin\big(\bigcup^\infty_{n=1}\partial(I_n\times J_n)\big)\cup\big(
\bigcap^\infty_{m=1}C_m\big).
\]
Then, there exists an $m_0\in\N$ such that $x\notin C_{m_0}$. Since $\bigcup\limits^{n_{m_0}-1}_{n=1}\partial(I_n\times J_n)$ is compact and clearly $x\notin\bigcup\limits^{n_{m_0}-1}_{n=1}\partial(I_n\times J_n)$, there exists a $\de(x)>0$ such that for every interval $A\times B\subseteq\R^2$ containing $x$, with $diam(A\times B)<\de(x)$, we have
\begin{eqnarray}
(A\times B)\cap\big(\bigcup^{n_{m_0}-1}_{n=1}(I_n\times J_n)\big)=\emptyset. \label{eq7}
\end{eqnarray}
We choose an $\el_0\in\N$ such that $n_{s_{\el_0}}>n_{m_0}$ and w.l.o.g. we take\linebreak $\de(x)<2^{s_{\el_0}}\cdot w_{n_{s_{\el_0}+1}-1}$.

Taking $t\in(0,\de(x))$, there exists some $\el_\ast\in\N$ with $\el_\ast>\el_0$ such that
\begin{eqnarray}
2^{s_{\el_\ast+1}}\cdot w_{n_{s_{\el_\ast+1}+1}-1}\le t<2^{s_{\el_\ast}}\cdot w_{n_{s_{\el_\ast}+1}-1}.\label{eq8}
\end{eqnarray}

For $A\times B$ with $diam(A\times B)<t$, since $x\in A\times B$ and $x\notin C_{m_0}$, by (\ref{eq5}), (\ref{eq6}) and [3, Proposition 3.6] we have that
\begin{eqnarray}
(A\times B)\cap\big(\bigcup^{n_{s_{\el_\ast}+1}-1}_{i=n_{m_0}}(I_i\times J_i)\big)=\emptyset.  \label{eq9}
\end{eqnarray}
By (\ref{eq7}) and (\ref{eq9}) we get
\begin{eqnarray}
(A\times B)\cap\big(\bigcup^{n_{s_{\el_\ast}+1}-1}_{i=1}(I_i\times J_i)\big)=\emptyset.  \label{eq10}
\end{eqnarray}

Also, since $x\notin C_{m_0}$ and $s_{\el_\ast}>m_0$ (so $n_{s_{\el_\ast}}>n_{m_0}$), from the definition of $C_{m_0}$ we deduce that
\[
x\notin2^s-dil\big(\bigcup^{n_{s+1}-1}_{i=n_s}(I_i\times J_i)\big) \qquad \text{for} \ \ s\ge s_{\el_\ast+1}.
\]
So, by [3, Proposition 3.4], we have
\begin{eqnarray}
\frac{\big|(A\times B)\cap\big(\bigcup\limits^{n_{s+1}-1}_{i=n_s}(I_i\times J_i)\big)\big|}{|A\times B|}<\frac{2}{2^s} \qquad \text{for} \ \ s\ge s_{\el_\ast+1}.  \label{eq11}
\end{eqnarray}
By (\ref{eq10}) and (\ref{eq11}) we get
\[
\frac{\big|(A\times B)\cap\big(\bigcup\limits^\infty_{i=1}(I_i\times J_i)\big)\big|}{|A\times B|}<\sum^\infty_{s=s_{\el_\ast+1}}\frac{2}{2^s}.
\]
Consequently
\[
\frac{|(A\times B)\cap K|}{|A\times B|}>1-\sum^\infty_{s=s_{\el_\ast+1}}\frac{2}{2^s}
\]
for $x\in A\times B$ with $diam(A\times B)<t$.\\
Thus, in view of (\ref{eq8}), we obtain the conclusion of the theorem. \qs
\end{Proof}
\section{An elementary estimate for the auxiliary function}\label{sec4}
\noindent

The role of the following proposition is crucial in this section. Here is used the fact that the index $a\{w^2_n:n\in\N\}$ is strictly positive.
\begin{proposition}\label{prop7}
We have
\[
\underset{n}{\lim\inf}\frac{\log n}{|\log w_n|}>0.
\]
\end{proposition}
\begin{Proof}
We recall that $a\{w^2_n:n\in\N\}=\underset{n}{\lim\inf}a_n>0$, where $n\Big(\dfrac{r_n}{n}\Big)^{a_n}=1$ and $r_n=\ssum^\infty_{m=n}w_m^2$. Consequently
\[
r_n=\bigg(\frac{1}{n}\bigg)^{(\frac{1}{a_n}-1)}
\]
and since $\underset{n}{\lim\inf}a_n>0$
\begin{eqnarray}
r_n\ge\bigg(\frac{1}{n}\bigg)^\mi \qquad \text{finally for every $n\in\N$ and some $\mi>0$}.  \label{eq12}
\end{eqnarray}
Next, since $e_{BT}\{w^2_n:n\in\N\}<1$ (see (\ref{eq2})), we have
\begin{eqnarray}
r_n\le(w^2_n)^{1-\de} \qquad \text{for some $e_{BT}\{w^2_n:n\in\N\}<\de<1$}.  \label{eq13}
\end{eqnarray}
By (\ref{eq12}) and (\ref{eq13}), we have finally for every $n\in\N$
\[
\bigg(\frac{1}{n}\bigg)^\mi\le r_n\le w^{2(1-\de)}_n
\]
equivalently
\[
\bigg(\frac{1}{n}\bigg)^{\mi/2(1-\de)}<w_n
\]
so
\[
\frac{\mi}{2(1-\de)}\cdot\log\bigg(\frac{1}{n}\bigg)<\log w_n
\]
and since $\log w_n<0$ (note that $w_n\ra0$)
\[
\frac{\mi}{2(1-\de)}\cdot\log n>|\log w_n|
\]
thus
\[
\frac{\log n}{|\log w_n|}>\frac{2(1-\de)}{\mi}>0. \text{\qs}
\]
\end{Proof}\bigskip

We introduce the following notation:\\
For the amounts $a(t)>0$, $\bi(t)>0$ (with $t>0$), we write
\[
a(t)\cong\bi(t) \qquad \text{for} \ \  t\ra0 \ \  \text{or} \ \  t\ra\infty
\]
iff
\[
0\lvertneqq \underset{t}{\lim\inf}\bigg(\frac{a(t)}{\bi(t)}\bigg)\le\underset{t}{\lim\sup}
\bigg(\frac{a(t)}{\bi(t)}\bigg)<+\infty \qquad \text{for} \ \ t\ra0 \ \ \text{or} \ \ t\ra\infty.
\]
\begin{lemma}\label{lem8}
Under the above notation we have that
\[
\log n\cong|\log w_n| \qquad \text{for} \ \ n\ra\infty.
\]
\end{lemma}
\begin{Proof}
Immediate from Proposition \ref{prop7} and the equation
\[
e_{BT}\{w_n:n\in\N\}=2\cdot e_{BT}\{w^2_n:n\in\N\}\bigg(=2\cdot\underset{n}{\lim\sup}\frac{\log n}{|\log w^2_n|}\bigg),
\]
since by [4, pp. 274-275] and [7, p. 292 \& THEOREM, p. 35]
\[
e_{BT}\{w_n:n\in\N\}=\underset{n}{\lim\sup}\frac{\log n}{|\log w_n|}.  \text{\qs}
\]
\end{Proof}
\begin{lemma}\label{lem9}

\[
\log n\cong\log(n+1) \qquad \text{for} \ \ n\ra\infty.
\]
\end{lemma}
\begin{Proof}
Obvious. \qs
\end{Proof}
\begin{lemma}\label{lem10}
The following holds
\[
|\log(w_{n_{s_{\el+1}+1}-1})|\cong s_{\el+1}\cdot\log(s_{\el+1})
\qquad{for} \ \ \el\ra\infty.
\]
\end{lemma}
\begin{Proof}
Since $n_s:=s^s$, for $s\in\N$, by Lemmata \ref{lem8} and \ref{lem9} we have
\begin{align*}
&|\log(w_{n_{s_{\el+1}+1}-1})|\cong\log(n_{s_{\el+1}+1}-1)\cong\log(n_{s_{\el+1}+1}) \\
&=\log\Big((s_{\el+1}+1)^{s_{\el+1}+1}\Big)=(s_{\el+1}+1)\cdot
\log(s_{\el+1}+1) \\
&\cong s_{\el+1}\cdot\log(s_{\el+1}).  \text{\qs}
\end{align*}
\end{Proof}

In view of the definition of the auxiliary function, we define\linebreak $f:(0,2^{s_1}\cdot w_{n_{s_1+1}-1})\ra\R$, as follows:
\[
f(t)=\sum^\infty_{k=s_{\el+1}}\frac{2}{2^k}=\frac{1}{2^{s_{\el+1}}}\cdot
\bigg(\frac{2}{1-\Big(\dfrac{1}{2}\Big)}
\bigg)=4\cdot\bigg(\frac{1}{2}\bigg)^{s_{\el+1}} \]
for $2^{s_{\el+1}}\cdot w_{n_{s_{\el+1}+1}-1}\le t<2^{s_\el}\cdot w_{n_{s_\el+1}-1}$ and $\el\in\N$.
\begin{lemma}\label{lem11}
\[
f(2^{s_{\el+1}}\cdot w_{n_{s_{\el+1}+1}-1})\cdot|\log(2^{s_{\el+1}}\cdot w_{n_{s_{\el+1}+1}-1})|\ra0
\qquad for \ \ \el\ra\infty
\]
(note that $2^{s_{\el+1}}\cdot w_{n_{s_{\el+1}+1}-1}\ra0$ \ \ for $\el\ra\infty$).
\end{lemma}
\begin{Proof}
We have
\begin{eqnarray}
|\log(2^{s_{\el+1}}\cdot w_{n_{s_{\el+1}+1}-1})|\le s_{\el+1}\cdot\log 2+|\log
(w_{n_{s_{\el+1}+1}-1})| \label{eq14}
\end{eqnarray}
and
\begin{eqnarray}
f(2^{s_{\el+1}}\cdot w_{n_{s_{\el+1}+1}-1})=4\cdot\bigg(\frac{1}{2}\bigg)^{s_{\el+1}}
\qquad \text{for} \ \ \el\in\N \label{eq15}
\end{eqnarray}
By (\ref{eq14}), (\ref{eq15}) and Lemma \ref{lem10}, we have the conclusion of the
lemma. \qs
\end{Proof}
\begin{proposition}\label{prop12}
It holds that
\[
f(t)=o\bigg(\frac{1}{|\log t|}\bigg) \qquad \text{for} \ \ t\ra0.
\]
\end{proposition}
\begin{Proof}
This follows immediately from Lemma \ref{lem11} and the fact that $\dfrac{1}{|\log t|}$ is
concave nearby zero, while $f$ is a step function. \qs
\end{Proof}

In conclusion we have
\begin{theorem}\label{thm13}
Let $K\subseteq\R^2$ be a compact set of positive Lebesgue measure, $I=(\al,b)
\times(c,d)$ such that $K\subseteq I$ and $K=I\big\backslash\dis\bigcup_{n\in\N}
(I_n\times J_n)$, where $I_n\times J_n$, $n\in\N$ are disjoint cubes,
$I_n=J_n=w_n$ with $\{w_n:n\in\N\}$ a non-increasing sequence, $a\{w^2_n:n\in\N\}\gneqq0$ and
$e_{BT}\{w^2_n:n\in\N\}\lneqq1$.
Then, for $h$ associated to the decomposition of
$K=I\big\backslash\dis\bigcup_{n\in\N}(I_n\times J_n)$ we have
\[
\inf_{x\in A\times B\atop diam(A\times B)<t}\frac{|(A\times B)\cap K|}{|A\times B|}\ge1-o
\bigg(\frac{1}{|\log t|}\bigg)
\]
for a.e. $x\in K$ and for sufficiently small $t$.
\end{theorem}
\begin{Proof}
This is a corollary of Theorem \ref{thm6}, Proposition \ref{prop7} and the
fact that
\[
h(t)=1-f(t).  \text{\qs}
\]
\end{Proof}

\end{document}